\documentclass[12pt,a4paper]{article}
\usepackage[colorlinks=true,linkcolor=black, citecolor=blue, urlcolor=blue]{hyperref}
\usepackage{epsfig}
\usepackage{graphics}
\usepackage{amsfonts}
\usepackage{amsmath}
\usepackage{latexsym,mathrsfs}

\newtheorem{theorem}{Theorem}[section]
\newtheorem{lemma}{Lemma}[section]
\newtheorem{proposition}{Proposition}[section]
\newtheorem{corollary}{Corollary}[section]

\newtheorem{example}{Example}[section]

\newtheorem{remark}{Remark}[section]

\DeclareMathOperator{\supp}{supp} 
 \numberwithin{equation}{section}

\title{\textbf{On the projections of the multifractal packing dimension for $q>1$}}
\author{Bilel Selmi}

\begin{document}
\maketitle
\begin{abstract}
The aim of this article is to study the behavior of the multifractal
packing function $B_\mu(q)$ under projections in Euclidean space for
$q>1$. We show that $B_\mu(q)$ is preserved under almost every
orthogonal projection. As an application, we study the multifractal
analysis of the projections of a measure. In particular,  we obtain
general results for the multifractal analysis of the orthogonal
projections on $m$-dimensional linear subspaces of a measure $\mu$
satisfying the multifractal formalism.
\end{abstract}

\noindent{\bf MSC-2010:} {28A20, 28A80.}


\section{Introduction and statement of the results}
The notion of singularity exponents or spectrum and generalized
dimensions are the major components of the multifractal analysis.
They were introduced with a view of characterizing the geometry of
measure and to be linked with the multifractal spectrum which is the
map which affects the Hausdorff or packing dimension of the
iso-H\"{o}lder set
$$
E_\mu(\alpha)=\left\{x\in \supp\mu;\;\; \lim_{r\to
0}\frac{\log\big(\mu B(x,r)\big)}{\log r}=\alpha \right\}
$$
for a given $\alpha\geq 0$ and $\supp\mu$ is the topological support
of probability measure $\mu$ on $\mathbb{R}^n$, $B(x,r)$ is the
closed ball of center $x$ and radius $r$. It unifies the
multifractal spectra to the multifractal packing function $B_\mu(q)$
via the Legendre transform \cite{BBH, OL}, i.e.,
$$
\dim_H\Big(E_\mu(\alpha)\Big)=\inf_{q\in\mathbb{R}}\Big\{q\alpha+B_{\mu}(q)\Big\}.
$$
There has been a great interest in understanding the fractal
dimensions of projections of the iso-H\"{o}lder sets and measures.
Recently, the projectional behavior of dimensions and multifractal
spectra of sets and measures have generated a large interest in the
mathematical literature \cite{BF, DB, FM1, FJ, JJ, J, J1, SP, SP1}.
The first significant work in this area was the result of Marstrand
\cite{M} who showed who proved a well-known theorem according to
which the Hausdorff dimension of a planar set is preserved under
orthogonal projections. This result was later generalized to higher
dimensions by Kaufman \cite{K}, Mattila \cite{MP} and Hu and Taylor
\cite{HT} and  they obtain similar results for the Hausdorff
dimension of a measure.

Let us mention that Falconer and Mattila \cite{FM} and Falconer and
Howroyd \cite{FH} have proved that the packing dimension of the
projected set or measure will be the same for almost all
projections. However, despite these substantial advances for fractal
sets, only very little is known about the multifractal structure of
projections of measures, except a paper by O'Neil  \cite{O}. Later,
in \cite{BB} Barral and Bhouri, studied the multifractal analysis of
the orthogonal projections on $m$-dimensional linear subspaces of
singular measures on $\mathbb{R}^n$ satisfying the multifractal
formalism. The result of O'Neil was later generalized by Selmi et
al. in \cite{DB, BD1, BD3, SBB1, SB2}.

\bigskip
O'Neil \cite{O} has compared the generalized Hausdorff and packing
dimensions of a set $E$ of $\mathbb{R}^n$ with respect to a measure
$\mu$ with those of their projections onto $m$-dimensional
subspaces. More specifically, given a compactly supported Borel
probability measure $\mu$ on $\mathbb{R}^n$ and $q\in \mathbb{R}$,
let $B_\mu(q)$ the multifractal packing function of $\supp\mu$. Then
we have $B_{\mu_V}(q)\leq B_\mu(q)$ for all $q\leq 1$ and all
$m$-dimensional linear subspaces $V$. Then, what can be said about
the multifractal packing function and its projection onto a lower
dimensional linear subspace for $q>1$? The goal of this work is
giving an answer to this question. We are interested in knowing
whether or not this property is preserved after orthogonal
projections on $\gamma_{n,m}$-almost every linear $m$-dimensional
subspaces for $q>1$, where  $\gamma_{n,m}$ is the uniform measure on
$G_{n,m}$, the set of linear $m$-dimensional subspaces of
$\mathbb{R}^n$ endowed with its natural structure of a compact
metric space (see \cite{MP1}).

\bigskip
In the present paper we pursue those kinds of studies and we
consider the  multifractal formalism developed in \cite{O}. The aims
of this study are twofold. First, the behavior of the  packing
dimensions $B_\mu(q)$ under projection. In particular, we show that
$B_\mu(q)$ is preserved under $\gamma_{n,m}$-almost every orthogonal
projection for $q>1$. We have treated an unsolved case by O'Neil
which is $q > 1$ and the result that we have obtained is optimal.
Secondly, to investigate a relationship between the multifractal
spectrum and its projection onto a lower dimensional linear
subspace. We also obtain general results for the multifractal
analysis of the orthogonal projections on $m$-dimensional linear
subspaces of a measure $\mu$ satisfying the multifractal formalism.

\section{Preliminaries}
 We start by recalling the multifractal formalism
introduced by O'Neil in \cite{O}. Let $\mu$ be a compactly supported
probability measure on $\mathbb{R}^n$. For $q, s\in\mathbb{R}$, $E
\subseteq\supp\mu$ and $\delta>0$, we define the multifractal
packing pre-measure,

$$\overline{{\mathscr P}}^{q,s}_{\mu,\delta}(E) =\displaystyle
\sup \left\{\sum_i \mu\left(B\left(x_i,\frac{r_i}3\right)\right)^q
r_i^s\right\},$$ where the supremum is taken over all
$\delta$-packings of $E$,
$$\overline{{\mathscr P}}^{q,s}_{\mu}(E)
=\displaystyle\inf_{\delta>0}\overline{{\mathscr
P}}^{q,s}_{\mu,\delta}(E).$$

 \bigskip
The function $\overline{{\mathscr P}}^{q,s}_{\mu}$ is increasing but
not $\sigma$-subadditive. That is the reason why O'Neil introduced
the modification of the multifractal packing measure ${\mathscr
P}^{q,s}_{\mu}$:
$$
{\mathscr P}^{q,s}_{\mu}(E) = \inf_{E \subseteq \bigcup_{i}E_i}
\sum_i \overline{\mathscr P}^{q,s}_{\mu}(E_i).
 $$
In a similar way, we define the Hausdorff measure,
$${{\mathscr
H}}^{{q},s}_{{\mu},\delta}(E) = \displaystyle\inf \left\{\sum_i
\mu\Big(B(x_i,3r_i)\Big)^q ~r_i^s;\;\;\Big(B(x_i,
r_i)\Big)_i\text{is a $\delta$-covering of}\; E\right\},
 $$
and
 $$
{{\mathscr H}}^{{q},s}_{{\mu}}(E) =
\displaystyle\sup_{\delta>0}{{\mathscr H}}^{q,s}_{{\mu},\delta}(E).
 $$

The functions ${\mathscr P}^{q,s}_{\mu}$ and ${\mathscr
H}^{q,s}_{\mu}$ are metric outer measures and thus measures on the
family of Borel subsets of $\mathbb{R}^n$. An important feature of
the pre-packing, packing and Hausdorff measure is that ${\mathscr
P}^{q,s}_{\mu}\leq{\overline{\mathscr P}}^{q,s}_{\mu}$ and there is
a constant $c$ depending also on the dimension of the ambient space,
such that ${\mathscr H}^{q,s}_{\mu}\leq c~{\mathscr P}^{q,s}_{\mu}$
(see \cite{O}).

\smallskip\smallskip
The functions $\overline{{\mathscr P}}^{q,s}_{\mu}$,  ${\mathscr
P}^{q,s}_{\mu}$ and ${\mathscr H}^{q,s}_{\mu}$ assign, in the usual
way, a dimension to each subset $E$ of $\supp\mu$. They are
respectively denoted by $\Lambda_{\mu}^q(E)$, $B_{\mu}^q(E)$ and
$b_{\mu}^q(E)$.

\par\noindent\begin{enumerate}\item  There exists a unique number
$\Lambda_{\mu}^q(E)\in[-\infty,+\infty]$ such that
 $$
\overline{{\mathscr P}}^{q,s}_{\mu}(E)=\left\{\begin{matrix}
\infty &\text{if}& s < \Lambda_{\mu}^q(E),\\ \\
0 & \text{if}&  \Lambda_{\mu}^q(E) < s.
\end{matrix}\right.
 $$
\item  There exists a unique number
$B_{\mu}^q(E)\in[-\infty,+\infty]$ such that
 $$
{\mathscr P}^{q,s}_{\mu}(E)=\left\{\begin{matrix}
\infty &\text{if}& s <B_{\mu}^q(E),\\
 \\
 0 & \text{if}&  B_{\mu}^q(E) < s.\end{matrix}\right.
 $$
 \item  There exists a unique number
$b_{\mu}^q(E)\in[-\infty,+\infty]$ such that
 $$
{\mathscr H}^{q,s}_{\mu}(E)=\left\{\begin{matrix}
\infty &\text{if}& s <b_{\mu}^q(E),\\
 \\
 0 & \text{if}&  b_{\mu}^q(E) < s.\end{matrix}\right.
 $$
\end{enumerate}
We note that for all $q\in\mathbb{R}$
$$
b_{\mu}^q(\emptyset)=B_{\mu}^q(\emptyset)=\Lambda_{\mu}^q(\emptyset)=-\infty,
$$
and if $\mu(E)=0$, then
$$
b_{\mu}^q(E)=B_{\mu}^q(E)=\Lambda_{\mu}^q(E)=-\infty\quad\text{for}\quad
q>0.
$$
Next, we define the separator functions $\Lambda_{\mu}$, $B_{\mu}$
and $b_{\mu}$ : $\mathbb{R}\rightarrow [-\infty, +\infty]$ by,
\begin{center}
$\Lambda_{\mu}(q)=\Lambda_{\mu}^q(\supp\mu)$,
$B_{\mu}(q)=B_{\mu}^q(\supp\mu)$ and
$b_{\mu}(q)=b_{\mu}^q(\supp\mu).$
\end{center}
It is well known that the functions $\Lambda_{\mu}$, $B_{\mu}$ and
$b_{\mu}$ are decreasing. The functions $\Lambda_{\mu}$, $B_{\mu}$
convex and satisfying $b_{\mu}\leq B_{\mu}\leq \Lambda_{\mu}.$

\begin{proposition}\cite{O}
Let $\mu$ be compactly supported probability measure on
$\mathbb{R}^n$. Then, we have
\begin{enumerate}
\item For $q<1$, $0\leq b_{\mu}(q)\leq B_{\mu}(q)\leq
\Lambda_{\mu}(q).$
\item $b_{\mu}(1)=B_{\mu}(1)=
\Lambda_{\mu}(1)=0.$
\item For $q>1$, $b_{\mu}(q)\leq B_{\mu}(q)\leq
\Lambda_{\mu}(q)\leq 0.$
\end{enumerate}
\end{proposition}
\begin{remark}
The multifractal Hausdorff and packing measures introduced by O'Neil
are different from those developed by Olsen \cite{OL}. Although,
when $\mu$ satisfies a doubling condition, the multifractal measures
are equivalent.
\end{remark}
\section{Main result}

Let $\mu$ be a compactly supported probability measure on
$\mathbb{R}^n$ and $q\in \mathbb{R}$. In the following, we require
an alternative characterization of the generalized upper
$L^q$-spectrum of $\mu$  in terms of a potential obtained by
convolving $\mu$ with a certain kernel. For this purpose let us
introduce some notations. For $1\leq s\leq n$ and $r>0$ we define
 $$
\begin{array}{llll}\label{j}
\phi_r^s: & \mathbb{R}^n & \longrightarrow & \mathbb{R} \\
& x & \longmapsto & \min\Big\{1,\: r^s|x|^{-s}\Big\},
\end{array}
$$
and
\begin{equation*}
\label{r}\mu\ast\phi_r^s(x)=\int\min\Big\{1,\:
r^s|x-y|^{-s}\Big\}d\mu(y).\end{equation*} Let $E$ be a compact
subset of $\supp\mu$. For $1\leq s\leq n$ and $q>1$, write
$$
N_{\mu,r}^{q,s}(E)=
\int_E\Big(\mu\ast\phi_{r/3}^s(x)\Big)^{q-1}d\mu(x),
 $$
and
$$
\overline{\tau}_{\mu}^{q,s}(E)=\displaystyle\limsup_{r\to0}
\frac{\log N_{\mu,r}^{q,s}(E)}{-\log r}\quad\text{and}\quad
\underline{\tau}_{\mu}^{q,s}(E)=\displaystyle\liminf_{r\to0}
\frac{\log N_{\mu,r}^{q,s}(E)}{-\log r}.
$$

The definition of these dimensions is, frankly, messy, indirect and
unappealing. In an attempt to make the concept more attractive, we
present here an alternative approach to the dimension
$\overline{\tau}_{\mu}^{q,s}$ and his application to projections in
terms of a potential obtained by convolving $\mu$ with a certain
kernel. For $E$ a compact subset of $\supp\mu$ we can try to
decompose $E$ into a countable number of pieces $E_1, E_2, . . .$ in
such a way that the largest piece has as small a dimension as
possible. The present approach was first used by Falconer in
\cite[Section 3.3]{Falconer} and further developed by O'Neil in
\cite[Proposition 2.4]{O}. This idea leads to the following modified
dimension in terms of the convolutions:
\begin{eqnarray*}
\mathfrak{T}_\mu^{q,s}(E)&=&\inf\left\{\sup_{1\leq
i<\infty}\overline{\tau}_{\mu}^{q,s}(E_i),\;\; E\subset\bigcup_i
E_i\;\;\text{with each $E_i$ compact }\right\}
\end{eqnarray*}
and
$$
\mathfrak{T}_\mu^{s}(q)=\mathfrak{T}_\mu^{q,s}(\supp\mu)\quad\text{for
all}\quad s\geq 1.
$$

 Let $m$ be an integer with $0<m\leq n$ and $G_{n,m}$ the
Grassmannian manifold of all $m$-dimensional linear subspaces of
$\mathbb{R}^n$. Denote by $\gamma_{n,m}$ the invariant Haar measure
on $G_{n,m}$ such that $\gamma_{n,m}(G_{n,m})=1$. For $V\in
G_{n,m}$, we define the projection map $\pi_V:
\mathbb{R}^n\longrightarrow V$ as the usual orthogonal projection
onto $V$. Then, the set $\{\pi_V,\; V \in G_{n,m}\}$ is compact in
the space of all linear maps from $\mathbb{R}^n$ to $\mathbb{R}^m$
and the identification of $V$ with $\pi_V$ induces a compact
topology for $G_{n,m}$. Also, for a Borel probability measure $\mu$
with compact support $supp\mu \subset\mathbb{R}^n$ and for $V\in
G_{n,m}$, we denote by $\mu_V$, the projection of $\mu$ onto $V$,
i.e.,
 $$
\mu_V(A)=\mu(\pi_V^{-1}(A))\quad \forall A\subseteq V.
 $$

Since $\mu$ is compactly supported and $supp\mu_V=\pi_V(supp\mu)$
for all $V\in G_{n,m}$, then, for any continuous function $f:
V\longrightarrow\mathbb {R}$, we have
 $$
\displaystyle\int_V fd\mu_V=\int f(\pi_V(x))d\mu(x),
 $$
whenever these integrals exist. Then for all $V\in G_{n,m}$, $x\in
\mathbb{R}^n$ and $0<r<1$, we have
\begin{equation*}
\label{r}\mu\ast\phi_r^m(x)=\int\mu_V(B(x_V,r))dV=\int\min\Big\{1,\:
r^m|x-y|^{-m}\Big\}d\mu(y).\end{equation*}

\bigskip
In \cite{O}, O'Neil has compared the generalized Hausdorff and
packing dimensions of a set $E$ of $\mathbb{R}^n$ with respect to a
measure $\mu$ with those of their projections onto $m$-dimensional
subspaces. More specifically, he proved the following result:
\begin{theorem}
Let $\mu$ be a compactly supported probability measure on
$\mathbb{R}^n$ and $E\subseteq\supp\mu$. For $q\leq 1$ and all $V\in
G_{n,m}$, we have
$$
B_{\mu_V}^{q}(\pi_V(E))\leq B_{\mu}^{q}(E).
$$
\end{theorem}

\bigskip
In this paper, we show that $B_\mu^q(E)$ is preserved under
$\gamma_{n,m}$-almost every orthogonal projection for $q>1$.  We
have treated an unsolved case by O'Neil which is $q > 1$ and the
result that we have obtained is optimal. More precisely, we have the
following result.
\begin{theorem} \label{TH1} Let $E$ be a compact subset of $\supp\mu$ and $q>1$.
\begin{enumerate}
\item If $1<q\leq2$, one has
$$
B_{\mu_V}^{q}(\pi_V(E))=\mathfrak{T}_{\mu}^{q,m}(E)=\max\Big(m(1-q),
B_{\mu}^{q}(E)\Big),\;\text{for $\gamma_{n,m}$-almost every}\;\;
V\in G_{n,m}.
$$
 \item If $q>2$ and $({E}_i)_i$  is a cover
of $E$ by a countable collection of compact sets is such that
$\Lambda_{\mu}^{q}(E_i)\geq -m$ for all $i$, then
 $$
B_{\mu_V}^{q}(\pi_V(E))=\mathfrak{T}_{\mu}^{q,m}(E)=B_{\mu}^{q}(E),\;\text{for
$\gamma_{n,m}$-almost every}\;\; V\in G_{n,m}.
$$
\end{enumerate}
\end{theorem}
\newpage
As a consequence we have, the following corollary
\begin{corollary}\label{C1} Let $q>1$.
\begin{enumerate}
\item If $1<q\leq2$, one has
$$
B_{\mu_V}(q)=\mathfrak{T}_{\mu}^{m}(q)=\max\Big(m(1-q),
B_{\mu}(q)\Big),\;\text{for $\gamma_{n,m}$-almost every}\;\; V\in
G_{n,m}.
$$
 \item If $q>2$ and  $({E}_i)_i$  is a cover
of $\supp\mu$ by a countable collection of compact sets is such that
$\Lambda_{\mu}^{q}(E_i)\geq -m$ for all $i$, then
 $$
B_{\mu_V}^{q}(q)=\mathfrak{T}_{\mu}^{m}(q)=B_{\mu}(q),\;\text{for
$\gamma_{n,m}$-almost every}\;\; V\in G_{n,m}.
$$
\end{enumerate}
\end{corollary}
\begin{remark} The hypothesis $\Lambda_{\mu}^{q}(E_i\cap E)\geq -m$ for all
$i$ implies that $\Lambda_{\mu}^{q}(E)\geq -m$. Nevertheless, we
don't know if the weaker condition $\Lambda_{\mu}^{q}(E)\geq -m$ is
sufficient to obtain the conclusion of Theorem \ref{TH1}.
\end{remark}
\section{Proof of the main result}

\subsection{Preliminary results}
We present the tools, as well as the intermediate results, which
will be used in the proof of our main result. Let $\mu$ be a
compactly supported probability measure on $\mathbb{R}^n$ and $q\in
\mathbb{R}$. We define the upper and lower $L^q$- spectrum of a
measure $\mu$.
 For a subset $E \subset
\supp\mu$, write
$$
N_{\mu,r}^q(E)= \sup\left\{ \displaystyle \sum_i
\mu\left(B\left(x_i,\frac r3\right)\right)^q;\;\;
\Big(B(x_i,r)\Big)_i \quad\text{is a packing of}\;\; E\right\}.
 $$
The upper respectively lower $L^q$- spectrum
$\overline{\tau}_{\mu}^q$ and $\underline{\tau}_{\mu}^q$ of $E$ is
defined by
$$
\overline{\tau}_{\mu}^q(E)=\displaystyle\limsup_{r\to0} \frac{\log
N_{\mu,r}^q(E)}{-\log
r}\quad\text{and}\quad\underline{\tau}_{\mu}^q(E)=\displaystyle\liminf_{r\to0}
\frac{\log N_{\mu,r}^q(E)}{-\log r}.
$$
By convention, if $r\in (0,1)$: $\frac{\log 0}{-\log r}=-\infty$.

 \bigskip
\noindent The following proposition is a consequence of the
multifractal formalism developed in \cite{O}.
\begin{proposition}\label{ah} Let $E$ be a compact subset of $\supp\mu$ and $q\in
\mathbb{R}$. One has
\begin{eqnarray*}
B_\mu^q(E)=\inf\left\{\sup_{1\leq
i<\infty}\overline{\tau}_{\mu}^{q}(E_i),\;\; E\subset\bigcup_i
E_i\;\;\text{with each $E_i$ compact}\;\right\}.
\end{eqnarray*}
\end{proposition}
Proposition \ref{ah} is a consequence from the following lemmas.
\begin{lemma}\label{L1} Let $E$ be a subset of $\supp\mu$ and $q\in
\mathbb{R}$. Then we have
$$
B_\mu^q(E)=\inf\left\{\sup_{1\leq
i<\infty}\overline{\tau}_{\mu}^{q}(E_i),\;\; E\subset\bigcup_i
E_i\right\}=\inf\left\{\sup_{1\leq
i<\infty}\Lambda_{\mu}^{q}(E_i),\;\; E\subset\bigcup_i E_i\right\}.
$$
\end{lemma}
\noindent{\bf Proof.} The Lemma is Proposition 2.4 of \cite{O}.
\begin{lemma}\noindent
\begin{enumerate}
\item For $q,s\in \mathbb{R}$, we have $\overline{{\mathscr P}}^{q,s}_{\mu}(E)=\overline{{\mathscr
P}}^{q,s}_{\mu}(\overline{E})$ and
$\Lambda_{\mu}^{q}(E)=\Lambda_{\mu}^{q}(\overline{E})$ for all
$E\subset \supp\mu$.
\item Let $E$ be a compact subset of $\supp\mu$ and $q\in
\mathbb{R}$. If
$$
\Lambda_{\mu}^{q}(E\cap V)=\Lambda_{\mu}^{q}({E})\quad\text{for all
open sets}\;\; V\;\;\text{with}\;\;  E\cap V \neq \emptyset,
$$
then,
$$
B_\mu^q(E)=\Lambda_{\mu}^{q}({E}).
$$
\end{enumerate}
\end{lemma}
\noindent{\bf Proof.} The first part is  Lemma 5.4.1 in \cite{O2}.
We will prove the second part. Let $E\subset \cup_i E_i$. Since
$E\subset \cup_i \overline{E}_i$, Baire's category theorem implies
that there exists an integer $k\in \mathbb{N}$ and an open set $V$
such that $\emptyset\neq E\cap V\subset  \overline{E}_k.$ Hence,
$$
\sup_i\Lambda_{\mu}^{q}(\overline{E}_i)\geq\Lambda_{\mu}^{q}(\overline{E}_k)\geq\Lambda_{\mu}^{q}(E\cap
V)=\Lambda_{\mu}^{q}({E}).
$$
Since the covering $(E_i)_i$ of $E$ was arbitrary, the previous
lemma now implies that
\begin{eqnarray*}
B_\mu^q(E)&=&\inf\left\{\sup_{1\leq
i<\infty}\Lambda_{\mu}^{q}(E_i),\;\; E\subset\bigcup_i
E_i\right\}\\&=&\inf\left\{\sup_{1\leq
i<\infty}\Lambda_{\mu}^{q}(\overline{E_i}),\;\; E\subset\bigcup_i
E_i\right\}\\&\geq&\Lambda_{\mu}^{q}({E}).
\end{eqnarray*}

The following straightforward estimates concern the behavior of
$\mu\ast\phi_r^n(x)$ as $r\to0$.
\begin{lemma}\cite{FO}\label{l1}
Let $1\leq m \leq n$ and $\mu$ be a compactly supported probability
measure on $\mathbb{R}^n$. For all $x\in \mathbb{R}^n$, we have
$$
c r^m\leq\mu\ast\phi_{r}^m(x)
$$
for all sufficiently small $r$, where $c>0$ is independent of $r$.
\end{lemma}
\begin{lemma}\cite{FO}\label{l2} Let $\mu$ be a compactly supported probability measure on
$\mathbb{R}^n$.
\begin{enumerate}
\item  For all $x\in \mathbb{R}^n$ and $r>0$
$$ \mu (B(x,r))\leq\mu\ast\phi_{r}^n(x).$$
\item Let $\varepsilon>0$. We have that for $\mu$-almost all $x$
 $$
r^{-\varepsilon} \mu (B(x,r))\geq\mu\ast\phi_r^n(x),$$ if $r$ is
sufficiently small.
\end{enumerate}
\end{lemma}

We use the properties of $\mu\ast\phi_r^n(x)$ to have a relationship
between the kernels and projected measures.
\begin{lemma}\cite{FO}\label{l3} Let $1\leq m \leq n$, $\mu$ be a compactly supported probability measure on
$\mathbb{R}^n$, $\varepsilon>0$ and $r$ is sufficiently small.
\begin{enumerate}
\item For all $V\in G_{n,m}$ and for $\mu$-almost all
$x\in\mathbb{R}^n$
 $$
{r}^\varepsilon\mu\ast\phi_{r}^m(x)\leq\mu_V(B(x_V, r)).
 $$

\item For $\gamma_{n,m}$-almost all $V\in
G_{n,m}$ and all $x\in\mathbb{R}^n$
 $$
{ r}^\varepsilon \mu\ast\phi_{r}^m(x)\geq\mu_V(B(x_V,r)).
 $$
\end{enumerate}
\end{lemma}

\smallskip\smallskip

The next result is essentially a restatement of \cite[Proposition
4.2]{BB} and \cite[Proposition 5.1]{B} (see also \cite[Lemma 2.6
(a)]{FO} and \cite{SS}). We provide a proof for the reader's
convenience.
\begin{proposition}\label{P1}
Let $E$ be a compact subset of $\supp\mu$. For $q>1$, we have
$$
\underline{\tau}_{\mu}^{q}(E)=\liminf_{r\to0}\frac{1}{-\log
r}\log\int_E\Big(\mu(B(x,r/3))\Big)^{q-1}d\mu(x)
$$
and
$$
\overline{\tau}_{\mu}^{q}(E)=\limsup_{r\to0}\frac{1}{-\log
r}\log\int_E\Big(\mu(B(x,r/3))\Big)^{q-1}d\mu(x).
$$
\end{proposition}
\noindent{\bf Proof.} Let $r>0$ and $\Big(B(x_i,r/3)\Big)_i$ be a
family of disjoint balls centered on $E$. Since $q>1$, we have
\begin{eqnarray*}
\int_E\Big(\mu(B(x,r/3))\Big)^{q-1}d\mu(x)&\geq&\int_{\bigcup_iB(x_i,r/9)}
\Big(\mu(B(x,r/3))\Big)^{q-1}d\mu(x)\\
&\geq&\sum_i \Big(\mu\left(B\left(x_i, r/9\right)\right)\Big)^q.
\end{eqnarray*}

On the other hand, for every $r>0$ we can apply Besicovitch's
covering theorem to $\Big(B(x,r/3)\Big)_{x\in E}$, to get a positive
integer $\xi(n)$ depending on $n$ only, as well as
$\mathcal{B}_1=\Big(B(x_{1,j},r/3)\Big)_j$,...,$\mathcal{B}_{\xi(n)}=\Big(B(x_{\xi(n),j},r/3)\Big)_j,$
$\xi(n)$ families of disjoint balls of radius $r/3$, such that
$E\subset \bigcup_{i=1}^{\xi(n)}\bigcup_jB(x_{i,j},r/3)$ and
\begin{eqnarray*}
\int_E\Big(\mu(B(x,r/9))\Big)^{q-1}d\mu(x)&\leq&\sum_{i=}^{\xi(n)}
\sum_j\int_{B(x_{i,j},r/3)}\Big(\mu(B(x,r/9))\Big)^{q-1}d\mu(x)\\
&\leq&\sum_{i=}^{\xi(n)} \sum_j\Big(\mu(B(x_{i,j},r/3))\Big)^{q}.
\end{eqnarray*}
Taking the logarithms and letting $r\to0$ yields the result.
\begin{remark}
Let $E$ be a compact subset of $\supp\mu$ and $q>1$. It is clear
that from Proposition \ref{P1} and Proposition 2.8 in \cite{FO},
$$
\overline{\tau}_{\mu}^{q}(E)=\overline{\tau}_{\mu}^{q,n}(E)
\quad\text{and}\quad\underline{\tau}_{\mu}^{q}(E)=\underline{\tau}_{\mu}^{q,n}(E).
$$
\end{remark}
With these definitions we have the following  corollaries.
\begin{corollary}
Let $E$ be a compact subset of $\supp\mu$. For $q>1$ and $1\leq
m\leq n$, we have
\begin{eqnarray*}
\underline{\tau}_{\mu}^{q,m}(E)&=&\inf\left\{s\geq0;\;\;
\limsup_{r\to0}r^s\int_E\Big(\mu\ast\phi_{r/3}^m(x)\Big)^{q-1}d\mu(x)=0\right\}
\end{eqnarray*}
and
\begin{eqnarray*}
\overline{\tau}_{\mu}^{q,m}(E)&=&\inf\left\{s\geq0;\;\;
\liminf_{r\to0}r^s\int_E\Big(\mu\ast\phi_{r/3}^m(x)\Big)^{q-1}d\mu(x)=0\right\}.
\end{eqnarray*}
\end{corollary}

\begin{corollary}
Let $E$ be a compact subset of $\supp\mu$. For $q>1$, we have
\begin{eqnarray*}
\underline{\tau}_{\mu}^{q}(E)&=&\inf\left\{s\geq0;\;\;
\limsup_{r\to0}r^s\int_E\Big(\mu(B(x,r/3))\Big)^{q-1}d\mu(x)=0\right\}
\end{eqnarray*}
and
\begin{eqnarray*}
\overline{\tau}_{\mu}^{q}(E)&=&\inf\left\{s\geq0;\;\;
\liminf_{r\to0}r^s\int_E\Big(\mu(B(x,r/3))\Big)^{q-1}d\mu(x)=0\right\}.
\end{eqnarray*}
\end{corollary}
\begin{proposition}\label{P2} Let $E$ be a compact subset of $\supp\mu$.
For $q>1$, we have
\begin{equation*}
\overline{\tau}_{\mu}^{q,m}(E)=\max\Big(m(1-q),
\overline{\tau}_{\mu}^{q}(E)\Big).
\end{equation*}
\end{proposition}
\noindent{\bf Proof.}
 Recalling that, from Lemma \ref{l2}
 $$
\mu\ast\phi_{ r/3}^m(x)\geq\mu(B(x, r/3)),
 $$
it will be clear that for $q>1$ we have
$$
\int_E\Big(\mu\ast\phi_{r/3}^m(x)\Big)^{q-1}d\mu(x)\geq\int_E\Big(\mu(B(x,
r/3))\Big)^{q-1}d\mu(x).
 $$
Hence, $$\underline{\tau}_{\mu}^{q,m}(E)\geq
\underline{\tau}_{\mu}^{q}(E)\quad\text{and}\quad\overline{\tau}_{\mu}^{q,m}(E)\geq
\overline{\tau}_{\mu}^{q}(E).$$ Without loss of generality, we may
assume that $\mu(E)>0$. From Lemma \ref{l1}, we get
\begin{eqnarray*}
\int_E\Big(\mu\ast\phi_{r/3}^n(x)\Big)^{q-1}d\mu(x) &\geq& c_1
r^{m(q-1)} \mu(E).
\end{eqnarray*}
Therefore, we obtain
$$
\underline{\tau}_{\mu}^{q,m}(E)\geq m(1-q).
$$

Let $t>\overline{\tau}_{\mu}^{q}(E)$.  Suppose that $\supp\mu$ have
diameter $h$. Then
 $$
\displaystyle\int_E\Big(\mu(B(x,r/3))\Big)^{q-1}d\mu(x)\leq
c_1r^{-t},\quad \forall r\leq 6h,
 $$
where $c_1$ is independent of $r$, and
 $$
\displaystyle\int_E\Big(\mu(B(x,r/3))\Big)^{q-1}d\mu(x)=1,\quad\forall
r\geq 6h.
 $$
For $\varepsilon>0$, $p=q-1$ and $r$ is small enough, by using
Proposition 2.5 in \cite{FO}, we obtain
\begin{eqnarray*}
\int_E\Big(\mu\ast\phi_{r/3}^m(x)\Big)^{p}d\mu(x) &\leq& B\;
(r/3)^{mp-\varepsilon}\int_{r/3}^{+\infty}u^{-mq-1}\int_E\Big(\mu(B(x,u))\Big)^pd\mu(x)du\\
&=& B\;
(r/3)^{mp-\varepsilon}\int_{r/3}^{6h}u^{-mq-1}\int_E\Big(\mu(B(x,u))\Big)^pd\mu(x)du\\
&+&B\;
(r/3)^{mp-\varepsilon}\int_{6h}^{+\infty}u^{-mp-1}\int_E\Big(\mu(B(x,u))\Big)^pd\mu(x)du\\
&\leq&C_1(r/3)^{mp-\varepsilon}\int_{r/3}^{6h}u^{-mp-1-t}du\\&+&C_2(r/3)^{mp-\varepsilon}
\int_{6h}^{+\infty} u^{-mp-1}du \\ \\&\leq&\left\{
\begin{array}{ll}
 \bigskip
C_3~~(r/3)^{-t-\varepsilon}   \quad\quad\quad\;   si & \hbox{$t>-mp$,} \\
\\
C_4~~(r/3)^{mp-\varepsilon}\quad\quad\quad si& t\leq -mp,
\end{array}\right.
\end{eqnarray*}
where $B$ and $C_i$ $(i=1,..., 4)$ are independent of $r$. This
gives that
 $$
\overline{\tau}_{\mu}^{q,m}(E)\leq\max(-mp, t),\quad \text{for
all}\quad t>\overline{\tau}_{\mu}^{q}(E).
 $$
Finally, we obtain
 $$
\overline{\tau}_{\mu}^{q,m}(E)\leq\max\Big(m(1-q),
\overline{\tau}_{\mu}^{q}(E)\Big).
 $$

\bigskip
The following results present alternative expressions of the
$L^q$-spectrum in terms of the convolutions as well as general
relations between the $L^q$-spectrum of a measure and that of its
orthogonal projections.
\begin{theorem} \label{TH2} Let $E$ be a compact subset of $\supp\mu$. Then, we have
\begin{enumerate}
\item for all $q>1$ and $V\in G_{n,m},$
 $$
\underline{\tau}_{\mu_V}^{q}(\pi_V(E))\geq\underline{\tau}_{\mu}^{q,m}(E)\quad\text{and}
\quad\overline{\tau}_{\mu_V}^{q}(\pi_V(E))\geq\overline{\tau}_{\mu}^{q,m}(E).
 $$

\item For all $1<q\leq2$ and $\gamma_{n,m}$-almost every  $V\in
G_{n,m},$
 $$
\overline{\tau}_{\mu_V}^{q}(\pi_V(E))=\overline{\tau}_{\mu}^{q,m}(E)=\max\Big(m(1-q),
\overline{\tau}_{\mu}^{q}(E)\Big)
 $$
and
 $$
\underline{\tau}_{\mu_V}^{q}(\pi_V(E))=\underline{\tau}_{\mu}^{q,m}(E).
 $$
 \item For all $q>2$ and $\gamma_{n,m}$-almost every  $V\in
G_{n,m},$
\begin{enumerate}
 \bigskip
\item If $-m\leq\overline{\tau}_{\mu}^{q}(E)$ then
$\overline{\tau}_{\mu_V}^{q}(\pi_V(E))=\overline{\tau}_{\mu}^{q,m}(E)=
\overline{\tau}_{\mu}^{q}(E).$
 \bigskip
\item $
\underline{\tau}_{\mu_V}^{q}(\pi_V(E))=\max\Big(
m(1-q),\underline{\tau}_{\mu}^{q,m}(E)\Big).$
\end{enumerate}
\end{enumerate}
\end{theorem}
\begin{remark} The assertion (2) is essentially a restatement
of the main result of Hunt et al. in \cite{HK} and  Falconer et al.
in \cite[Theorem 3.9]{FO}. The assertion (3) extends the result of
Hunt and Kaloshin (of Falconer and O'Neil) to the case $q> 2$
untreated in their work.
\end{remark}

\noindent{\bf Proof.}~ The first and second parts follows from
Proposition \ref{P2} and the following lemma which is a consequence
of Lemma \ref{l3}.
\begin{lemma}\label{15}
Let $E$ be a compact subset of $\supp\mu$. Then, we have
\begin{enumerate}
\item for all $q>1$ and $V\in G_{n,m},$
\begin{equation*}
\displaystyle\limsup_{r\longrightarrow 0}\frac{1} {-\log r}\left[
\log\left(\displaystyle\frac{\displaystyle\int_{E}\Big(\mu\ast\phi_{r/3}^n(x)\Big)^{q-1}d\mu(x)}
{\displaystyle\int_{\pi_V(E)}\Big(\mu_V(B(x_V,r/3))\Big)^{q-1}d\mu_V(x_V)}
\right)\right]\leq 0,
\end{equation*}
\item for $1<q\leq 2$ and $\gamma_{n,m}$-almost every $V\in G_{n,m},$
\begin{equation*}
\displaystyle\lim_{r\rightarrow 0}\frac{1} {-\log r}\left[
\log\left(\displaystyle\frac{\displaystyle\int_E\Big(\mu\ast\phi_{r/3}^n(x)\Big)^{q-1}d\mu(x)
}{\displaystyle\int_{\pi_V(E)}\Big(\mu_V(B(x_V,r/3))\Big)^{q-1}d\mu_V(x_V)}
\right)\right]=0.
\end{equation*}
\end{enumerate}
\end{lemma}
See \cite[Theorem 2.1]{FB},  \cite[Theorem 4.1]{B} and \cite{SB} for
the key ideas needed to prove the third part of Theorem \ref{TH2}.

\subsection{ Proof of Theorem \ref{TH1}} Let us prove our main
theorem. Let $q>1$.
\begin{enumerate}
\item If $s>\mathfrak{T}_{\mu}^{q,m}(E)$ we may cover $E$ by a countable
collection of sets $E_i$, which we may take to be compact, such that
$\overline{\tau}_{\mu}^{q,m}(E_i)<s$. By using Theorem \ref{TH2}
(2.), we have $\overline{\tau}_{\mu_V}^{q}(\pi_V(E_i))\leq s$ for
$\gamma_{n,m}$-almost every  $V\in G_{n,m}.$ Proposition \ref{ah}
implies that ${B}_{\mu_V}^{q}(\pi_V(E))\leq s$ for
$\gamma_{n,m}$-almost every  $V\in G_{n,m}$ and so,
${B}_{\mu_V}^{q}(\pi_V(E))\leq \mathfrak{T}_{\mu}^{q,m}(E)$ for
$\gamma_{n,m}$-almost every  $V\in G_{n,m}.$

Now, if $s<\mathfrak{T}_{\mu}^{q,m}(E)$. Fix $V\in G_{n,m}$ and let
$(\widetilde{E}_i)_i$ be a cover of the compact set $\pi_V(E)$ by a
countable collection of compact sets. Put for each $i$,
$E_i=E\cap\pi_V^{-1}(\widetilde{E}_i),$ then
$\sup_i\overline{\tau}_{\mu}^{q,m}(E_i)>s$. By using Theorem
\ref{TH2} (1.), we have
$\sup_i\overline{\tau}_{\mu_V}^{q}(\pi_V(E_i))\geq s$ and
$\sup_i\overline{\tau}_{\mu_V}^{q}(\widetilde{E}_i)\geq s$, this
implies that ${B}_{\mu_V}^{q}(\pi_V(E))\geq s$. Therefore, we obtain
${B}_{\mu_V}^{q}(\pi_V(E))\geq \mathfrak{T}_{\mu}^{q,m}(E)$.

Thus, the part concerning the equality between $\max\Big(m(1-q),
B_{\mu}^{q}(E)\Big)$ and $\mathfrak{T}_{\mu}^{q,m}(E)$ is a
consequence of Proposition \ref{P2}.

\item Let
$({E}_i)_i$  be a cover of $E$ by a countable collection of compact
sets is such that
$\overline{\tau}_{\mu}^{q}(E_i)=\Lambda_{\mu}^{q}(E_i)\geq -m$ for
all $i$. Then, by using Lemma \ref{L1} and Proposition \ref{P2}, we
have $\mathfrak{T}_{\mu}^{q,m}(E)=B_{\mu}^{q}(E)$.

Now, if $s>\mathfrak{T}_{\mu}^{q,m}(E)$ we may cover $E$ by a
countable collection of sets $E_i$, which we may take to be compact,
such that $\overline{\tau}_{\mu}^{q,m}(E_i)<s$. By using Theorem
\ref{TH2} (2.) and since
$-m\leq\overline{\tau}_{\mu}^{q}(E_i)=\Lambda_{\mu}^{q}(E_i)$ for
all $i$, we have $\overline{\tau}_{\mu_V}^{q}(\pi_V(E_i))\leq s$ for
$\gamma_{n,m}$-almost every  $V\in G_{n,m}.$ Proposition \ref{ah}
implies that ${B}_{\mu_V}^{q}(\pi_V(E))\leq s$ for
$\gamma_{n,m}$-almost every  $V\in G_{n,m}$ and so,
${B}_{\mu_V}^{q}(\pi_V(E))\leq \mathfrak{T}_{\mu}^{q,m}(E)$ for
$\gamma_{n,m}$-almost every  $V\in G_{n,m}.$ In similar way, we
prove ${B}_{\mu_V}^{q}(\pi_V(E))\geq \mathfrak{T}_{\mu}^{q,m}(E)$
for all $V\in G_{n,m}.$
\end{enumerate}

\bigskip\bigskip
We can improve substantially the O'Neil's result \cite[Corollary
5.12]{O} in the following example:
\begin{example}
Fix $0<m\leq n$ and let $\mu$ be a self-similar measure on
$\mathbb{R}^n$ with support equal to $K$ such that $\dim_P(K)=s\leq
m$. Let $q\geq0$ and $({E}_i)_i$  be a cover of $E$ by a countable
collection of compact sets is such that $\Lambda_{\mu}^{q}(E_i)\geq
-m$ for all $i$. By using Corollary \ref{C1} and Corollary 5.12 in
\cite{O}, we have for $\gamma_{n,m}$-almost every $V\in G_{n,m}$
$$
B_{{\mu}_V}({q})=b_{{\mu}_V}({q})= b_{{\mu}}({{q}})=B_{{\mu}}({q}).
 $$
\end{example}
\section{Application}
When  $\mu$ obeys the multifractal formalism over some interval, we
are interested in knowing whether or not this property is preserved
after orthogonal projections on $\gamma_{n,m}$-almost every linear
$m$-dimensional subspaces.

This section is devoted to study the behavior of projections of
measures obeying to the multifractal formalism. More precisely, we
prove that for $q>1$ if the multifractal formalism holds for $\mu$
at $\alpha=-B'_{\mu}(q)$, it holds for $\mu_V$ for
$\gamma_{n,m}$-almost every $V\in G_{n,m}$. Before detailing our
results let us recall the multifractal formalism introduced by
O'Neil. For $\alpha\geq 0$, let
$$
{E}_{\mu}(\alpha)=\left\{x\in\supp\mu;\;\;
\lim_{r\to0}\frac{\log\Big(\mu B(x,3r)\Big)}{\log r}=\alpha
\right\}.
$$
We mention that in the last decade there has been a great interest
for the multifractal analysis and positive results have been written
in various situations (see for example \cite{BBJ, BBH, OL, OL1}).

\bigskip
The function $B_\mu(q)$ is related to the multifractal spectrum of
the measure $\mu$. More precisely,
$f^*(\alpha)=\displaystyle\inf_\beta\big(\alpha\beta+f(\beta)\big)$
denotes the Legendre transform of the function $f,$ it has been
proved in \cite{BBJ, BBH, OL, OL1} a lower and upper bound estimate
of the singularity spectrum using the Legendre transform of the
function $B_\mu(q)$. The following theorem is a consequence of the
multifractal formalism developed in \cite{BBH}.
\begin{theorem}\label{th3}
Let $\mu$ be a compactly supported Borel probability measure on
$\mathbb{R}^n$ and $ q\in\mathbb{R}$. Suppose that
\begin{enumerate}
\item $\mathscr{H}_{\mu}^{q,B_{\mu}(q)}(\supp\mu)>0,$

\item $B_{\mu}$ is differentiable at $q$.
\end{enumerate}
Then,
\begin{eqnarray*}
\dim_H E_{\mu}\big(-B'_{\mu}(q)\big)=\dim_P
E_{\mu}\big(-B'_{\mu}(q)\big) =
B^*_{\mu}\big(-B'_{\mu}(q)\big)=b^*_{\mu}\big(-B'_{\mu}(q)\big).
\end{eqnarray*}
Here $\dim_H$ and $\dim_P$ denote, respectively, the Hausdorff and
the packing dimension, see \cite{OL} for precise definitions of
this.
\end{theorem}

The following proposition is established in \cite{O}.
\begin{proposition}\label{P}
Let $\mu$ be a compactly supported Borel probability measure on
$\mathbb{R}^n$. For $q\geq1$ and all $V\in G_{n,m}$, we have
 $$
b_{\mu_V}(q)\geq \max\Big(m(1-q), b_\mu(q)\Big).
 $$
\end{proposition}

 \bigskip
In the following, we study the validity of the multifractal
formalism under projection. More specifically, we obtain general
result for the multifractal analysis of the orthogonal projections
on $m$-dimensional linear subspaces of  measure $\mu$ satisfying the
multifractal formalism.
\begin{theorem}\label{TH3}
Let $\mu$ be a compactly supported Borel probability measure on
$\mathbb{R}^n$ and $ q>1$. Suppose that
 \bigskip
\\$(\mathsf H_1)$ $\mathscr{H}_{\mu}^{q,B_{\mu}(q)}(\supp\mu)>0,$
 \bigskip
\\$(\mathsf H_2)$ $B_{\mu}$ is differentiable at $q$,
 \bigskip
 \\$(\mathsf H_3)$ $({E}_i)_i$  be a cover
of $\supp\mu$ by a countable collection of compact sets is such that
$b_{\mu}^q(E_i\cap\supp\mu)\geq \max(-m, m(1-q))$ for all $i$.
 \bigskip
\\Then, for $\gamma_{n,m}$-almost every $V\in
G_{n,m}$,
\begin{eqnarray*}
\dim_P {E}_{\mu_V}\big(-B'_{\mu}(q)\big)&=& \dim_H
{E}_{\mu_V}\big(-B'_{\mu}(q)\big) =\dim_H
E_{\mu}\big(-B'_{\mu}(q)\big)\\&=& \dim_P
E_{\mu}\big(-B'_{\mu}(q)\big) = B^*_{\mu}\big(-B'_{\mu}(q)\big) =
b^*_{\mu}\big(-B'_{\mu}(q)\big).
\end{eqnarray*}
\end{theorem}
\begin{remark}
The results of Theorem \ref{TH3} hold if we replace the condition
$$\mathscr{H}_{\mu}^{q,B_{\mu}(q)}(\supp\mu)>0$$ by  the existence of
a nontrivial (Frostman) measure $\nu_q$ satisfying
$$
\nu_q(B(x,r))\leq  \mu(B(x,3r))^q ~r^{B_\mu(q)}
$$
where $x\in \supp\mu$ and $0 < r < 1$. \\For more details, the
reader can see \cite[Theorem 5.1]{O}.
\end{remark}
\noindent {\bf Proof.} By using Corollary \ref{C1}, Proposition
$\ref{P}$, $(\mathsf H_1)$ and $(\mathsf H_3)$ we have, for
$\gamma_{n,m}$-almost every $V\in G_{n,m}$,
\begin{equation}\label{qq}
b_{\mu}(q)=b_{\mu_V}(q)=B_{\mu_V}(q)=B_{\mu}(q).
\end{equation}
$(\mathsf H_1)$,  (\ref{qq}) and the proof of Lemma 3.2 in \cite{O}
ensure that, there exists a positive constant  $c$ such that
\begin{center}
$0<\mathscr{H}_{\mu}^{q,B_{\mu}(q)}(\supp\mu)\leq
c~\mathscr{H}_{\mu_V}^{q,B_{\mu_V}(q)}(\supp\mu_V),$\quad for
$\gamma_{n,m}$-almost every $V\in G_{n,m}$.
\end{center}
So, the hypothesis $(\mathsf H_2)$, Theorem \ref{th3} and the
equalities (\ref{qq}) imply that
\begin{eqnarray}\label{4}
\dim_H{E}_{\mu_V}\big(-B'_{\mu}(q)\big)\geq
-qB'_{\mu}(q)+B_{\mu}(q),
\quad\text{for}\;\;\gamma_{n,m}\text{-almost every} \;\;V\in
G_{n,m}.
\end{eqnarray}
Hence, the assumption (\ref{qq}) give that
\begin{eqnarray}\label{3}
\dim_P{E}_{\mu_V}\big(-B'_{\mu}(q)\big)&\leq&
B^*_{\mu_V}\big(-B'_{\mu}(q)\big)\nonumber\\
&=& B^*_{\mu}\big(-B'_{\mu}(q)\big),
\end{eqnarray}
for $\gamma_{n,m}$-almost every $V\in G_{n,m}$. Thus, the result is
a consequence from (\ref{4}) and (\ref{3}).

\section*{Acknowledgments} The author is greatly indebted to the referee for his/her
carefully reading the first submitted version of this paper and giving
elaborate comments and valuable suggestions on revision so that
the presentation can be greatly improved.

\bigskip\bigskip
\noindent {\bf Bilel SELMI}\\
\noindent Analysis, Probability and Fractals Laboratory LR18ES17\\
Faculty of Sciences of Monastir\\
Department of Mathematics\\
University of Monastir\\
5000-Monastir\\
Tunisia\\ \\
E-mail:\; bilel.selmi@fsm.rnu.tn
\end{document}